\documentclass[a4paper,11pt]{article}
\usepackage{ams}
\usepackage[T1]{fontenc}

\newtheorem{df}{Definition}[section]
\newtheorem{thm}[df]{Theorem}
\newtheorem{prop}[df]{Proposition}

\newtheorem{ex}[df]{Example}
\newtheorem{lem}[df]{Lemma}

\def\A{{\cal A}}

\def\boom{\quad\lower3pt\hbox{\vrule height1.1ex width .9ex depth -.2ex}
                    \vskip9pt}

\newcommand{\pf}{\noindent{\sc Proof.}\ }

\let\Bar=\overline

\let\Tilde=\widetilde

\def\til0{\widetilde 0}

\def\tila{\widetilde a}
\def\tilb{\widetilde b}

\def\tilq{\widetilde q}

\let\da=\partial
\let\isom=\cong
\let\sol=\bullet

\def\chigh{{\raise1.5pt\hbox{$\chi$}}}
\let\Ga=\Gamma
\let\Om=\Omega

\let\phi=\varphi

\def\plusH{\ \lower 5pt\hbox{${\buildrel {\textstyle +}
\over {\scriptscriptstyle H}}$}\ }
\def\minusH{\ \lower 5pt\hbox{${\buildrel {\textstyle -}
\over {\scriptscriptstyle H}}$}\ }
\def\timesH{\ \lower 4pt\hbox{${\buildrel {\textstyle .}
\over{\scriptscriptstyle H}}$}\ }

\def\plusV{\ \lower 5pt\hbox{${\buildrel {\textstyle +}
\over {\scriptscriptstyle V}}$}\ }
\def\minusV{\ \lower 5pt\hbox{${\buildrel {\textstyle -}
\over {\scriptscriptstyle V}}$}\ }
\def\timesV{\ \lower 4pt\hbox{${\buildrel {\textstyle .}
\over{\scriptscriptstyle V}}$}\ }

\def\dminus{\raise2pt\hbox{\vrule height1pt width 2ex}\hskip3pt}

\def\timesK{\ \lower 4pt\hbox{${\buildrel {\textstyle .}
\over{\scriptscriptstyle K}}$}\ }

\def\llangle{\langle\!\langle}
\def\rrangle{\rangle\!\rangle}

\def\LAvb{${\cal LA}$-vector bundle}
\def\LAgpd{${\cal LA}$-groupoid}

\def\CDO{\mathop{\rm CDO}}

\def\id{{\rm id}}

\def\ts#1{\stackrel{\vee}{#1}}

\def\pback#1{\mathbin{{{\lower1.2ex\hbox{$\times$}}\atop #1}}}
\def\sdp{\mathbin{\hbox{$\mapstochar\kern-.3333em\times$}}}

\def\gpd{\,\lower1pt\hbox{$\longrightarrow$}\hskip-.24in\raise2pt
             \hbox{$\longrightarrow$}\,}

\def\sgpd{\,\lower1pt\hbox{$\mlra$}\hskip-0.4in\raise2pt\hbox{$\mlra$}\,}

\def\vgpd{\Bigg\downarrow\!\!\Bigg\downarrow}

\def\vlra{\hbox{$\,-\!\!\!-\!\!\!-\!\!\!-\!\!\!-\!\!\!
-\!\!\!-\!\!\!-\!\!\!-\!\!\!-\!\!\!\longrightarrow\,$}}

\def\lgpd{\,\lower1pt\hbox{$\vlra$}\hskip-1.02in\raise2pt\hbox{$\vlra$}\,}

\def\mlra{\hbox{$\,-\!\!\!-\!\!\!\longrightarrow\,$}}

\def\hcompo#1#2#3{{\vcenter{\vbox{\hrule height.#2pt\hbox{\vrule width.#2pt
   height#1pt\kern#3pt\vrule width.#2pt\kern#3pt\vrule width.#2pt}
   \hrule height.#2pt}}}}

\def\vcompo#1#2#3{{\vcenter{\vbox{\hrule height.#2pt
   \hbox{\vrule width.#2pt height#3pt\kern#1pt\vrule width.#2pt}
   \hrule height.#2pt
   \hbox{\vrule width.#2pt height#3pt\kern#1pt\vrule width.#2pt}
   \hrule height.#2pt}}}}

\def\dsq{\mathop{\lower1pt\vbox{\hrule height.4pt \hbox
{\vrule width.4pt height.6em
\kern.6em \vrule width.4pt} \hrule height.4pt}}}

\def\dcomp{\mathop{\dsq\hskip-.88em\raise1pt\hbox{$\scriptstyle\nwarrow$}}}

\def\ssq{\vbox{\hrule height.4pt \hbox{\vrule width.4pt height.7in
\kern.7in \vrule width.4pt} \hrule height.4pt}}

\def\tsq{\mathop{\lower1pt\vbox{\hrule height.4pt \hbox
{\vrule width.4pt height.7em
\kern.7em \vrule width.4pt} \hrule height.4pt}}}

\def\m@th{\mathsurround=0pt}

\def\n@space{\nulldelimiterspace=0pt \m@th}

\def\Bbigg#1{{\hbox{$\left#1\vbox to 27.5pt{}\right.\n@space$}}}

\def\affel#1#2#3#4#5#6#7#8{\matrix{#7&\raise.7ex\hbox{$#1$}&#6\cr
                  \raise.5in\hbox{$#2$}&\sq&\raise.5in\hbox{$#3$}\cr
                          #8&#4&#5\cr}}

\begin{document}

\title{{\bf DOUBLE LIE ALGEBROIDS AND THE DOUBLE OF A LIE BIALGEBROID}
}

\author{K. C. H. Mackenzie\\
        School of Mathematics and Statistics\\
        University of Sheffield\\
        Sheffield, S3 7RH\\
        England\\
        {\sf K.Mackenzie@sheffield.ac.uk}}

\date{{\sf July 29, 1998}}

\maketitle

In this paper we define an abstract notion of double Lie algebroid,
and consider three major classes of examples.

Firstly, we verify that the double Lie algebroid of a double Lie groupoid,
and more generally an \LAgpd, as constructed in \cite{Mackenzie:Doubla2},
is an abstract double Lie algebroid. A large part of the necessary work for
this was done in \cite{Mackenzie:SDGDPG}.

Secondly we consider Lie bialgebroids. We show that the double cotangent
of a Lie bialgebroid is a double Lie algebroid and, further, that the
double cotangent of an {\em a priori} unrelated pair of Lie algebroid
structures on a vector bundle and its dual form a Lie bialgebroid if and
only if the double cotangent is a double Lie algebroid. We argue that this
is an appropriate form of the Manin triple result for Lie bialgebroids.

Thirdly we consider vacant double Lie algebroids and show that they
are equivalent to a matched pair structure on the two side Lie algebroids.

The paper begins with a preliminary study of the triple structures
associated with the tangent and cotangent of a double vector bundle.

Several of the results of this paper were announced without proof in
\cite{Mackenzie:DDEDLA}, which should be read as an introduction to this
paper. I am very grateful to Johannes Huebschmann, Yvette
Kosmann--Schwarzbach, Alan Weinstein and Ping Xu for conversations
on this material at various stages.

\newpage

\section{COTANGENT TRIPLE VECTOR BUNDLES}

A {\em double vector bundle} is a diagram as in Figure~\ref{fig:dvb}(a),
in which each side has a vector bundle structure, and the two structures
\begin{figure}[hbt]
\begin{picture}(395,100)
\put(20,40){$\matrix{       && \tilq_H           &&\cr
                            &{\cal A}&\mlra &A^V&\cr
                            &&&&\cr
                    \tilq_V &\Bigg\downarrow&&\Bigg\downarrow&q_V\cr
                            &&&&\cr
                            &A^H&\mlra &M&\cr
                            &&q_H&&\cr
                            &&&&\cr
                            &&\mbox{(a)}&&\cr
                            }$}
\put(270,40){$\matrix{&& \tilq^{(*)}_V &&\cr
        &{\cal A}^{*V}&\mlra   &K^*&\cr
        &&&&\cr
\tilq_{*V} &\Bigg\downarrow& &\Bigg\downarrow&q_{K^*} \cr
        &&&&\cr
        &A^H&\mlra &M&\cr
        &&q_H&&\cr
        &&&&\cr
        &&\mbox{(b)}&&\cr
        }$}
\end{picture}         \caption{\ \label{fig:dvb}}
\end{figure}
on ${\cal A}$ {\em commute} in the sense that the maps
defining each structure on ${\cal A}$ (the bundle projection, zero
section, addition and scalar multiplication) are morphisms with respect
to the other. This is precisely what is needed to ensure that when
four elements $\xi_1,\dots,\xi_4$ of ${\cal A}$ are such that the LHS of
$$
(\xi_1\plusH\xi_2)\plusV(\xi_3\plusH\xi_4) =
(\xi_1\plusV\xi_3)\plusH(\xi_2\plusH\xi_4)
$$
is defined, then the RHS is also, and they are equal. See
\cite{Pradines:DVB} or \cite[\S1]{Mackenzie:1992}.

The {\em core} of the double vector bundle ${\cal A}$ is the intersection
$K$ of the kernels of the two projections $\tilq_H$ and $\tilq_V$; the
vector bundle structures on ${\cal A}$ induce a common vector bundle
structure on $K$ with base $M$. The dual ${\cal A}^{*V}$ of the vertical
bundle structure on ${\cal A}$ has, in addition to its standard structure
on base $A^H$, a vector bundle structure on base $K^*$. The projection
is defined by
\begin{equation}                             \label{eq:unfproj}
\langle\tilq^{(*)}_V(\Phi), k\rangle =
      \langle\Phi, \til0^V_X \plusH \Bar k\rangle
\end{equation}
where $\Phi\colon \tilq_V^{-1}(X)\to\R,\ X\in A_m^H,$ and $k\in K_m$.
The addition $\plusH$ in ${\cal A}^{*V}\to K^*$ is defined by
\begin{equation}
\langle\Phi\plusH\Phi', \xi\plusH\xi'\rangle =
   \langle\Phi,\xi\rangle + \langle\Phi',\xi'\rangle
\end{equation}
and the zero above $\kappa\in K^*_m$ is $\til0^{(*V)}_\kappa$ defined by
$$
\langle\,\til0^{(*V)}_\kappa, \til0^H_x \plusV \Bar k\rangle =
\langle\kappa, k\rangle
$$
where $x\in A^V_m, k\in K_m.$ The scalar multiplication is defined in a
similar way. These two structures make ${\cal A}^{*V}$ a double vector
bundle as in Figure~\ref{fig:dvb}(b), the {\em vertical dual of ${\cal A}$}.
Its core is $(A^V)^*$: the core element $\Bar\psi$ corresponding to
$\psi\in(A^V_m)^*$ is
$$
\langle\Bar\psi, \til0^H_x \plusV \Bar k\rangle =
\langle\psi,x\rangle.
$$
See \cite{Pradines:1988} or \cite[\S3]{Mackenzie:SDGDPG}.

There is also a {\em horizontal dual} ${\cal A}^{*H}$ with sides $A^V$ and
$K^*$ and core $(A^H)^*$. The two duals are themselves dual, the pairing
being given by
\begin{equation}                       \label{eq:3duals}
\langle\Phi, \Psi\rangle = \langle\Psi, \xi\rangle
                            -  \langle\Phi, \xi\rangle
\end{equation}
where $\Phi\in {\cal A}^{*V},\ \Psi\in {\cal A}^{*H}$ have
$\tilq_V^{(*)}(\Phi) = \tilq_H^{(*)}(\Psi)$ and $\xi$ is any element
of ${\cal A}$ with $\tilq_V(\xi) = \tilq_{*V}(\Phi)$ and
$\tilq_H(\xi) = \tilq_{*H}(\Psi).$ See \cite[\S3]{Mackenzie:SDGDPG};
this pairing has also been found in \cite{KoniecznaU}.

This pairing could equally well be replaced by its negative. We regard
the choice of sign as an extra structure on ${\cal A}$ and write
$({\cal A};A^H,A^V;M)^+$ or $({\cal A};A^V,A^H;M)^-$ for the above
convention and $({\cal A};A^H,A^V;M)^-$ or $({\cal A};A^V,A^H;M)^+$
for the opposite.

For an ordinary vector bundle $(A,q,M)$ there is the tangent double vector
bundle of Figure~\ref{fig:TA}(a);
\begin{figure}[hbt]
\begin{picture}(395,100)
\put(20,40){$\matrix{&&T(q)  &&\cr
        &TA&\mlra&TM&\cr
        &&&&\cr
     p_A&\Bigg\downarrow&&\Bigg\downarrow&p\cr
        &&&&\cr
        &A&\mlra&M&\cr
        &&q&&\cr
        &&&&\cr
        &&\mbox{(a)}&&\cr
        }$}
\put(270,40){$\matrix{&& r    &&\cr
        &T^*A&\mlra&A^*&\cr
        &&&&\cr
     c_A&\Bigg\downarrow&&\Bigg\downarrow&q_*\cr
        &&&&\cr
        &A&\mlra&M&\cr
        &&q&&\cr
        &&&&\cr
        &&\mbox{(b)}&&\cr
        }$}
\end{picture}         \caption{\ \label{fig:TA}}
\end{figure}
for more detail on this see \cite{Besse} or \cite[\S5]{MackenzieX:1994}.
Its vertical dual is the cotangent double vector bundle of
Figure~\ref{fig:TA}(b). Its horizontal dual we denote
$(T^\sol A;A^*,TM;M)$; this is canonically isomorphic to
$(T(A^*);A^*,TM;M)$ under an isomorphism $I\colon T(A^*)\to T^\sol A$
given by
$$
\langle I({\cal X}),\xi \rangle =
\llangle{\cal X}, \xi \rrangle
$$
where ${\cal X}\in T(A^*), \xi \in TA$ have $T(q_*)({\cal X}) = T(q)(\xi)$,
and $\llangle\ ,\ \rrangle$ is the tangent pairing of $T(A^*)$ and $TA$
over $TM$. See \cite[\S5]{MackenzieX:1994}.

The structure of ${\cal A} = TA$ thus induces a pairing of $T^*A$ and
$T(A^*)$ over $A^*$ given by
$$
\langle\Phi,{\cal X}\rangle =
\llangle{\cal X},\xi\rrangle -
\langle\Phi, \xi\rangle
$$
where $\Phi\in T^*A, {\cal X}\in T(A^*)$ have $r(\Phi) = p_{A^*}({\cal X})$
and $\xi\in TA$ is chosen so that $T(q)(\xi) = T(q_*)({\cal X})$ and
$p_A(\xi) = c_A(\Phi)$. This pairing is nondegenerate by a general result
\cite[3.1]{Mackenzie:SDGDPG} so it defines an isomorphism of double vector
bundles $R\colon T^*A^*\to T^*A$ by the condition
$$
\langle R({\cal F}), {\cal X}\rangle = \langle{\cal F}, {\cal X}\rangle
$$
where the pairing on the RHS is the standard one of $T^*(A^*)$ and
$T(A^*)$ over $A^*$. This $R$ preserves the side bundles $A$ and $A^*$
but induces $-\id\colon T^*M\to T^*M$ as the map of cores. In summary we
now have the very useful equation
\begin{equation}                                   \label{eq:vue}
\langle{\cal F}, {\cal X}\rangle + \langle R({\cal F}), \xi\rangle
= \llangle{\cal X}, \xi\rrangle,
\end{equation}
for ${\cal F}\in T^*A^*, {\cal X}\in T(A^*), \xi\in TA$, where the pairings
are over $A^*, A$ and $TM$ respectively. See \cite[5.5]{MackenzieX:1994}.

Now return to the general double vector bundle in Figure~\ref{fig:dvb}(a),
denoting the core by $K$. Since ${\cal A}$ has two vector bundle
structures, it has two double cotangent bundles.
These fit together into a triple structure as the left and rear faces of
Figure~\ref{fig:cottrip}(a), the top face being essentially the cotangent
double of the two duals of ${\cal A}$.
\begin{figure}[htb]
\begin{picture}(350,200)
\put(0,150){$\matrix{&&      &\cr
                      &T^*{\cal A}&\mlra &{\cal A}^{*H}\cr
                      &&&\cr
                      &\Bigg\downarrow   & &\Bigg\downarrow     \cr
                      &&&\cr
                      &{\cal A}&\mlra &A^V\cr}$}

\put(40, 160){\vector(3,-4){40}}                
\put(110, 160){\vector(3,-4){40}}               
\put(40, 100){\vector(3,-4){40}}                
\put(120, 100){\vector(3,-4){30}}               

\put(75,70){$\matrix{&&      &\cr
                     &{\cal A}^{*V} &\mlra &K^*\cr
                     &&&\cr
                     &\Bigg\downarrow &&\Bigg\downarrow \cr
                     &&&\cr
                     &A^H &\mlra & M\cr}$}

\put(100,0){(a)}


\put(200,150){$\matrix{&&      &\cr
                      &T{\cal A}&\mlra &T(A^V)\cr
                      &&&\cr
                      &\Bigg\downarrow & &\Bigg\downarrow \cr
                      &&&\cr
                      &{\cal A}&\mlra &A^V \cr}$}

\put(240, 160){\vector(3,-4){40}}                
\put(315, 170){\vector(3,-4){40}}               
\put(240, 100){\vector(3,-4){40}}                
\put(320, 100){\vector(3,-4){30}}               

\put(265,70){$\matrix{&&      &\cr
                     &T(A^H) &\mlra &TM\cr
                     &&&\cr
                     &\Bigg\downarrow &&\Bigg\downarrow \cr
                     &&&\cr
                     &A^H &\mlra & M\cr}$}
\put(300,0){(b)}
\end{picture}\caption{\ \label{fig:cottrip}}
\end{figure}
This is, in a sense we will make precise elsewhere, the vertical dual
of the tangent prolongation of ${\cal A}$ in Figure~\ref{fig:cottrip}(b).
Five of the six faces in (a) are double vector bundles of types considered
already; it is only necessary to verify that the top face is a double vector
bundle.
(In all diagrams of this type, we take the oblique arrows to be coming
out of the page.)
The cores of the five faces are known, and we take the core of the
top face to be $T^*K$, in accordance with \cite[1.5]{Mackenzie:SDGDPG}.
Taking these cores in pairs, with the edges parallel to them, then gives three
double vector bundles: the left--right cores form $(T^*(A^H);A^H,(A^H)^*;M)$,
the back--front cores form $(T^*(A^V);A^V,(A^V)^*;M)$ and the top--down
cores form $(T^*K;K,K^*;M)$. Each of these {\em core double vector bundles}
has core $T^*M$.

We will also need to consider the cotangent triples of the two duals
of ${\cal A}$. Figure~\ref{fig:cofd}(a) is the cotangent triple of the
double vector bundle $({\cal A}^{*V};A^H,K^*;M)$; the ${}^\dagger$
denotes the dual over $K^*$.
\begin{figure}[htb]
\begin{picture}(350,200)
\put(0,150){$\matrix{&&     &&\cr
                     &T^*({\cal A}^{*V})&\mlra &({\cal A}^{*V})^\dagger & \isom {\cal A}^{*H}\cr
                     &&&&\cr
                     &\Bigg\downarrow   & &\Bigg\downarrow  &   \cr
                     &&&&\cr
                     &{\cal A}^{*V}&\mlra &K^*&\cr}$}

\put(40, 160){\vector(3,-4){40}}                
\put(120, 165){\vector(3,-4){40}}               
\put(40, 100){\vector(3,-4){40}}                
\put(120, 100){\vector(3,-4){30}}               

\put(75,70){$\matrix{&&      &\cr
                     &{\cal A} &\mlra &A^V\cr
                     &&&\cr
                     &\Bigg\downarrow &&\Bigg\downarrow \cr
                     &&&\cr
                     &A^H &\mlra & M\cr}$}

\put(100,0){(a)}


\put(200,150){$\matrix{&&      &\cr
                      &T({\cal A}^{*V})&\mlra &T(K^*)\cr
                      &&&\cr
                      &\Bigg\downarrow & &\Bigg\downarrow \cr
                      &&&\cr
                      &{\cal A}^{*V}&\mlra &K^* \cr}$}

\put(240, 160){\vector(3,-4){40}}                
\put(320, 165){\vector(3,-4){40}}               
\put(240, 100){\vector(3,-4){40}}                
\put(320, 100){\vector(3,-4){30}}               

\put(265,70){$\matrix{&&      &\cr
                     &T(A^H)&\mlra &TM\cr
                     &&&\cr
                     &\Bigg\downarrow &&\Bigg\downarrow \cr
                     &&&\cr
                     &A^H &\mlra & M\cr}$}
\put(300,0){(b)}
\end{picture}\caption{\ \label{fig:cofd}}
\end{figure}
We use the isomorphisms of double vector bundles
$$
Z_V\colon ({\cal A}^{*H})^\dagger \to {\cal A}^{*V}, \qquad
Z_H\colon ({\cal A}^{*V})^\dagger \to {\cal A}^{*H}
$$
induced by the pairing (\ref{eq:3duals}). Note that $Z_V$ preserves both
sides, $A^H$ and $K^*$, but induces $-\id\colon (A^V)^*\to (A^V)^*$ on the
cores, while $Z_H$ preserves $K^*$ and the core $(A^H)^*$, but induces
$-\id$ on the sides $A^V$; this reflects the fact that $Z_V = Z_H^\dagger$,
the dual over $K^*$; see \cite[3.6]{Mackenzie:SDGDPG}.

In the case of Figure~\ref{fig:TA}(a), we have $Z_V = R\circ I^\dagger$,
where the ${}^\dagger$ dual here is over $A^*$.

\section{ABSTRACT DOUBLE LIE ALGEBROIDS}

We come now to the definition of a double Lie algebroid. It will be
useful to have a name for a very special case.

\begin{df}
An {\em \LAvb}\ is a double vector bundle as in Figure~{\em \ref{fig:dvb}(a)}
together with Lie algebroid structures on a pair of parallel sides, such
that the structure maps of the other pair of vector bundle structures are
Lie algebroid morphisms.
\end{df}

For definiteness, take the Lie algebroid structures to be on
${\cal A}\to A^H$ and $A^V\to M$.

In the terminology of \cite[\S4]{Mackenzie:1992}, an \LAvb\ is an \LAgpd\
in which the groupoid structures are vector bundles (and in which the
scalar multiplication also preserves the Lie algebroid structures).
The core of an \LAgpd\ has a Lie algebroid structure induced from the
Lie algebroid structure on ${\cal A}$ \cite[\S5]{Mackenzie:1992}. Each
$k\in\Ga K$ induces $\Bar{k}\in\Ga_{A^H}{\cal A}$ defined by
$\Bar{k}(X) = k(m)\plusH\til0^V_X$ for $X\in A^H_m$ and the bracket on
$\Ga K$ is obtained by $\Bar{[k,\ell]} = [\Bar{k}, \Bar{\ell}]$.

\begin{lem}                                    \label{lem:K0}
The anchor and the bracket on the core of an \LAvb\ are zero.
\end{lem}

\pf
The anchor $\tila_V\colon{\cal A}\to T(A^H)$ is a morphism of double vector
bundles and therefore induces a core map $\da^H\colon K\to A^H$ which, by
\cite[\S5]{Mackenzie:1992}, is a Lie algebroid morphism. Since $A^H$ is
abelian, we have $a_K = a_H\circ\da^H = 0$.

Horizontal scalar multiplication by $t\neq 0$ defines a morphism
${\cal A}\to{\cal A}$ over $A^H\to A^H$ and therefore induces a map of
sections $t_H\colon \Ga_{A^H}{\cal A} \to \Ga_{A^H}{\cal A}$; the Lie
algebroid condition then ensures that
$[t_H(\xi), t_H(\eta)] = t_H([\xi,\eta])$
for all $\xi, \eta\in\Ga_{A^H}{\cal A}$. Now for $k\in\Ga K$,
$\Bar{tk} = t_H(\Bar{k})$ and so $\Bar{t[k,\ell]} = [t_H(\xi), t_H(\eta)]
= [\Bar{tk}, \Bar{t\ell}] = \Bar{t^2[k,\ell]}$. Therefore the bracket
must be zero.
\boom

We can also apply the calculus developed in \cite[\S3]{Mackenzie:SDGDPG}
for general \LAgpd s. Consider the Poisson structure on ${\cal A}^{*V}$.
Since it is linear over $A^H$, the Poisson anchor
$\pi^{\#V}\colon T^*({\cal A}^{*V})\to T({\cal A}^{*V})$ is a morphism
of double vector bundles for the left faces of Figure~\ref{fig:cofd},
with the corner map ${\cal A}\to T(A^H)$ being $\tila_V$ and core map
$-\tila_V^*\colon T^*(A^H)\to {\cal A}^{*V}.$

Now applying \cite[3.14]{Mackenzie:SDGDPG}, the Poisson structures on
${\cal A}^{*V}$ and $K^*$ make ${\cal A}^{*V}\to K^*$ a Poisson groupoid;
since the Poisson structure on $K^*$ is zero, this is a Poisson vector
bundle in the usual sense. Thus $\pi^{\#V}$ is also a morphism of
double vector bundles for the rear faces of Figure~\ref{fig:cofd}.
Denote the corner map $({\cal A}^{*V})^\dagger\to T(K^*)$ by $\chigh_V$;
since $\pi^{\#V}$ is skew--symmetric, the core map of the rear faces is
$-\chigh_V^*\colon T^*K^*\to {\cal A}^{*V}$. It then follows by a simple
argument (as in \cite[2.3]{Mackenzie:SDGDPG}) that $\pi^{\#V}$ is a
morphism of triple vector bundles.

We now turn to the general notion of double Lie algebroid. Again consider a
double vector bundle as in Figure~\ref{fig:dvb}(a). We now assume that there
are Lie algebroid structures on all four sides. The definition comprises
three conditions.

\subsection{Condition I}

{\em With respect to the two vertical Lie algebroids, ${\cal A}\to A^H$
and $A^V\to M$, the double vector bundle ${\cal A}$ is an \LAvb. Likewise,
with respect to the two horizontal Lie algebroids, ${\cal A}\to A^V$
and $A^H\to M$, the double vector bundle ${\cal A}$ is an \LAvb.}

\bigskip

Denote the four anchors by $\tila_V\colon {\cal A} \to T(A^H),\
\tila_H\colon {\cal A}\to T(A^V),\ a_V\colon A^V\to TM$ and
$a_H\colon A^H\to TM$. As usual we denote all four brackets by $[\ ,\ ]$;
the notation for elements will make clear which structure we are using.

The anchors thus give morphisms of double vector bundles
$$
(\tila_V;\id,a_V;\id)\colon ({\cal A};A^H,A^V;M)\to (T(A^H);A^H,TM;M),
$$
$$
(\tila_H;a_H,\id;\id)\colon ({\cal A};A^H,A^V;M)\to (T(A^V);TM,A^V;M)
$$
and so define morphisms of their cores; denote these by
$\da^H\colon K\to A^H$ and $\da^V\colon K\to A^V$.

Now return to $\pi^{\#V}$ and Figure~\ref{fig:cofd}.
Since the corner map ${\cal A}\to T(A^H)$ is $\tila_V$, the corner map
$A^V\to TM$ is $a_V$ (or $-a_V$ if $Z_H$ is incorporated) and the core
map for the front faces is $\da^H$. Likewise, since the core map of the
left faces is $-\tila^*_V$, the core map of the right faces must be
$-(\da^H)^*\colon (A^H)^*\to K^*$ (whether or not $Z_H$ is incorporated).
Lastly, the core map of the top faces is
$\pi^\#_V\colon T^*((A^V)^*)\to T((A^V)^*)$, the anchor for the Poisson
structure on $(A^V)^*$ dual to the given Lie algebroid structure on $A^V$.
(These observations are all special cases of \cite[\S3]{Mackenzie:SDGDPG}.)

Each of the maps of the core double vector bundles induces on
$T^*M\to (A^V)^*$ the map $-a_V^*$.

Similarly we can analyze $\pi^{\#H}\colon T^*({\cal A}^{*H})\to
T({\cal A}^{*H})$ as a morphism of triple vector bundles, and obtain
$\chigh_H\colon ({\cal A}^{*H})^\dagger\to T(K^*)$.

For Condition II, note first that it is automatic that $\tila_V$ is a
morphism of Lie algebroids over $A^H$ and that $a_V$ is a morphism
of Lie algebroids over $M$.

\subsection{Condition II}

{\em The anchors $\tila_V$ and $a_V$ form a morphism of Lie algebroids with
respect to the horizontal structure on ${\cal A}$ and the prolongation to
$TA^H\to TM$ of the structure on $A^H\to M$. Likewise, the anchors $\tila_H$
and $a_H$ form a morphism of Lie algebroids with respect to the vertical
structure on ${\cal A}$ and the prolongation to $TA^V\to TM$ of the structure
on $A^V\to M$.
}

\bigskip

By Condition I and the discussion preceding it,
the Poisson structure on ${\cal A}^{*H}\to K^*$ is
linear, and therefore induces a Lie algebroid structure on its dual
$({\cal A}^{*H})^\dagger$. We use $Z_V$ to transfer this to
${\cal A}^{*V}\to K^*$. Similarly the linear Poisson structure on
${\cal A}^{*V}\to K^*$ induces a Lie algebroid structure on
$({\cal A}^{*V})^\dagger\to K^*$.

\subsection{Condition III}

{\em With respect to these structures,
$({\cal A}^{*V}, ({\cal A}^{*V})^\dagger)$ is a Lie bialgebroid.
Further, $({\cal A}^{*V}; AH, K^*; M)$ is an \LAvb\ with respect to the
horizontal Lie algebroid structures and \break $({\cal A}^{*H}; K^*, AV ;M)$
is an \LAvb\ with respect to the vertical structures.
}

\bigskip

\begin{df}                                    \label{df:doubla}
A {\em double Lie algebroid} is a double vector bundle as in
Figure~{\em \ref{fig:dvb}(a)} equipped with Lie algebroid structures on all
four sides such that the above conditions {\em I, II, III} are satisfied.
\end{df}

The notion of Lie bialgebroid was defined in \cite{MackenzieX:1994} in
terms of the coboundary operators associated to $A$ and to $A^*$; a more
efficient and elegant reformulation was then given in
\cite{Kosmann-Schwarzbach:1995}. The definition most useful to us here
is quoted below in \ref{thm:6.2}. For the moment we only need the
following.

Suppose that $(E, E^*)$ is a Lie bialgebroid on base $P$ and denote the
anchors by $e$ and $e_*$. Then we take the Poisson structure on $P$ to
be $\pi^\#_P = e_*\circ e^*$; this is the opposite to \cite{MackenzieX:1994},
but the same as \cite{Kosmann-Schwarzbach:1995}. It follows that $e$ is a
Poisson map (to the tangent lift structure on $TP$) and $e_*$ is
anti--Poisson.

One expects the core of a double Lie algebroid to have a Lie algebroid
structure induced from those on ${\cal A}$. However, as \ref{lem:K0}
shows, the straightforward embedding of $K$ in terms of core sections
yields only the zero structure (see also \cite{Mackenzie:Doubla2}). Here
we obtain the correct structure in terms of its dual.

The anchor of $({\cal A}^{*V})^\dagger$ is $\chigh_V$, the appropriate
corner map of the Poisson anchor for ${\cal A}^{*V}$. On the other hand,
the anchor for ${\cal A}^{*V}$ itself is $\chigh_H\circ Z_V^{-1}.$ So
the Poisson anchor for $K^*$ is
$$
\pi^\#_{K^*} = \chigh_V\circ (Z_V^{-1})^\dagger \circ \chigh_H^\dagger.
$$
Now $Z_V^\dagger = Z_H$ has side map $-\id\colon A^V\to A^V$ and
$\chigh_V$ has side map $a_V$. The side map of $\chigh_H^\dagger$ is the
dual of the core map $-(\da^V)^*$ of $\chigh_H$. Thus the side map of
$\pi^\#_{K^*}$ is $a_V\circ\da^V$. One likewise checks that the core map
is $-(\da^H)^*\circ a_H^*$. This proves the first half of the following
result.

\begin{prop}
The anchor $a_K$ for the Lie algebroid structure on $K$ induced by the
Poisson structure on $K^*$ which arises from the Lie bialgebroid structure
on $({\cal A}^{*V}, ({\cal A}^{*V})^\dagger)$ is
$a_V\circ\da^V = a_H\circ\da^H$. The maps $\da^H\colon K\to A^H$ and
$\da^V\colon K\to A^V$ are Lie algebroid morphisms.
\end{prop}

\pf
Since $\chigh_V$ is the anchor for $({\cal A}^{*V})^\dagger$, it is
anti--Poisson into $T(K^*)$. Regarding $\chigh_V$ as a morphism of
the right faces in Figure~\ref{fig:cofd}, its core is $-(\da^H)^*$,
which is therefore anti--Poisson. So $\da^H$ is a morphism of Lie
algebroids.
\boom

The most fundamental example motivating \ref{df:doubla} is of course
that of the double Lie algebroid of a double Lie groupoid, as
constructed in \cite{Mackenzie:1992}, \cite{Mackenzie:Doubla2}. Most
of what is required to verify that the double Lie algebroid of a double
Lie groupoid does satisfy \ref{df:doubla} has been proved in
\cite{Mackenzie:SDGDPG}, but we recall the details briefly.

In order to proceed, we need to describe the notion of double Lie groupoid
in more detail (see \cite{Mackenzie:1992} and references given there).
A double Lie groupoid consists of a manifold $S$ equipped with two Lie
groupoid structures on bases $H$ and $V$, each of which is a Lie
groupoid over a common base $M$, such that the structure maps (source,
target, multiplication, identity, inversion) of each groupoid structure
on $S$ are morphisms with respect to the other; see Figure~\ref{fig:S}(a).
\begin{figure}[htb]
\begin{picture}(395,100)
\put(20,50){$\matrix{&&&&\cr
          &S&\sgpd &V&\cr
          &&&&\cr
          &\vgpd&&\vgpd&\cr
          &&&&\cr
          &H&\sgpd &M&\cr
          &&&&\cr
          &&\mbox{(a)}&&\cr}$}
\put(140,50){$\matrix{      &&   &&\cr
                            &A_VS    &\sgpd &AV&\cr
                            &&&&\cr
                            &\Bigg\downarrow&&\Bigg\downarrow&\cr
                            &&&&\cr
                            &H &\sgpd &M&\cr
                            &&&&\cr
                            &&\mbox{(b)}\cr}$}
\put(280,50){$\matrix{&&&&\cr
        &A^2S&\mlra   &AV&\cr
        &&&&\cr
        &\Bigg\downarrow& &\Bigg\downarrow&   \cr
        &&&&\cr
        &AH&\mlra &M&\cr
        &&&&\cr
        &&\mbox{(c)}\cr
        }$}
\end{picture}         \caption{\ \label{fig:S}}
\end{figure}
One should think of elements of $S$ as squares, the horizontal edges of
which come from $H$, the vertical edges from $V$, and the corner points
from $M$.

Consider a double Lie groupoid $(S;H,V;M)$ as in Figure~\ref{fig:S}(a).
Applying the Lie functor to the vertical structure $S\gpd H$ gives a Lie
algebroid $A_VS\to H$ which has also a groupoid structure over $AV$ obtained
by applying the Lie functor to the structure maps of $S\gpd V$; this is
the {\em vertical \LAgpd}\ of $S$ \cite[\S4]{Mackenzie:1992}, as in
Figure~\ref{fig:S}(b). The Lie algebroid of $A_VS\gpd AV$ is denoted $A^2S$;
there is a double vector bundle structure $(A^2S;AH,AV;M)$ obtained by
applying $A$ to the vector bundle structure of $A_VS\to H$
\cite{Mackenzie:Doubla2}; see Figure~\ref{fig:S}(c). Reversing the order of
these operations, one defines first the {\em horizontal \LAgpd}\
$(A_HS;AH,V;M)$ and then takes the Lie algebroid $A_2S = A(A_HS)$. The
canonical involution $J_S\colon T^2S\to T^2S$ then restricts to an
isomorphism of double vector bundles $j_S\colon A^2S\to A_2S$ and
allows the Lie algebroid structure on $A^2S\to AV$ to be transported to
$A_2S\to AV$. Thus $A_2S$ is a double vector bundle equipped with four
Lie algebroid structures; in \cite{Mackenzie:Doubla2} we called this
the {\em double Lie algebroid of $S$}. The core of both double vector bundles
$A_2S$ and $A^2S$ is $AC$, the Lie algebroid of the core groupoid $C\gpd M$
of $S$ \cite[1.6]{Mackenzie:Doubla2}.

Consider $A_2S$. The structure maps for the horizontal vector bundle
$A_2S\to AV$ are obtained by applying the Lie functor to the structure
maps of $A_HS\to V$ and are therefore Lie algebroid morphisms with
respect to the vertical Lie algebroid structures. The corresponding
statement is true for the vertical vector bundle $A^2S\to AH$ and this
is transported by $j_S$ to $A_2S\to AH$. Thus Condition~I holds.

Let $\ts{a}_V\colon A_2S\to TAH$ denote the anchor for the Lie algebroid of
$A_HS\gpd AH$. Then, as with any Lie groupoid, $\ts{a}_V = A(\ts{\chigh}_V)$
where $\ts{\chigh}_V\colon A_HS\to AH\times AH$ combines the target and
source of $A_HS\gpd AH$. It is easily checked that $\ts{\chigh}_V$ is a
morphism of \LAgpd s over $\chigh_V\colon V\to M\times M$ and $\id\colon
AH\to AH$, and so it follows, by using the methods of
\cite[\S1]{Mackenzie:Doubla2}, that $\ts{a}_V$ is a morphism of Lie algebroids
over $a_V$. Similarly one transports the result for the anchors
$A^2S\to TAV$ and $AH\to TM$. Thus Condition~II is satisfied.

Now consider the bialgebroid condition. The vertical dual
${\cal A}^{*V}$ is $A^*(A_HS)$ and in order to take the dual of this
over $A^*C$ we use the isomorphism ${j'}^V\colon A^*(A_HS)\to A(A^*_VS)$ of
\cite[(20)]{Mackenzie:SDGDPG}. This induces $({\cal A}^{*V})^\dagger
\isom A^*(A^*_VS)$.

Now the structure on ${\cal A}^{*V}$ itself comes from
$({\cal A}^{*H})^\dagger$. We have ${\cal A}^{*H} = A^\sol(A_HS)$ and
the isomorphism $I_H\colon A(A^*_HS)\to A^\sol(A_HS)$ associated to the
\LAgpd\ $A_HS$ in \cite[\S3]{Mackenzie:SDGDPG} (see also (\ref{eq:needI})
below) allows us to identify $({\cal A}^{*H})^\dagger$ with
$A^*(A^*_HS)$.

So $({\cal A}^{*V}, ({\cal A}^{*V})^\dagger)$ is effectively given by
$(A^*(A^*_HS), A^*(A^*_VS))$. Now use the isomorphism
${\cal D}_H\colon A^*(A^*_HS)\to A(A^*_VS)$ of \cite[3.9]{Mackenzie:SDGDPG}
and we have $(A(A^*_VS), A^*(A^*_VS))$; this is the Lie bialgebroid of
$A^*_VS\gpd A^*C$ which was proved to be a Poisson groupoid in
\cite[2.12]{Mackenzie:SDGDPG}. Notice that we started with $A_2S$,
defined in terms of the horizontal \LAgpd, and ended with the Lie
bialgebroid of the dual of the vertical \LAgpd.

To make this sketch into a proof, one must ensure that the various
isomorphisms preserve the Poisson structures involved. Rather than do
this, we prove a more general result.

Consider an \LAgpd\ as in Figure~\ref{fig:LAgpd}(a); that is, $\Om$ is
both a Lie algebroid over $G$ and a Lie groupoid over $A$, and each of the
groupoid structure maps is a Lie algebroid morphism; further, the map
$\Om\to A\pback{M} G$ defined by the source and the bundle projection, is a
surjective submersion.
\begin{figure}[htb]
\begin{picture}(395,100)
\put(20,40){$\matrix{      &&\tilq   &&\cr
                            &\Om     &\mlra &G &\cr
                            &&&&\cr
                            &\vgpd&&\vgpd&\cr
                            &&&&\cr
                            &A  &\mlra &M&\cr
                            &&q_A&&\cr
                            &&&&\cr
                            &&\mbox{(a)}&&\cr
                            }$}
\put(280,40){$\matrix{      &&A(\tilq) &&\cr
                            &A\Om     &\mlra &AG&\cr
                            &&&&\cr
                   \ts{q}   &\Bigg\downarrow&&\Bigg\downarrow&q_G\cr
                            &&&&\cr
                            &A  &\mlra &M&\cr
                            &&q_A&&\cr
                            &&&&\cr
                            &&\mbox{(b)}&&\cr
                            }$}
\end{picture}         \caption{\ \label{fig:LAgpd}}
\end{figure}
Applying the Lie functor vertically gives a double vector bundle
${\cal A} = A\Om$ as in Figure~\ref{fig:LAgpd}(b), with Lie algebroid
structures on the vertical sides. It is shown in \cite[\S1]{Mackenzie:Doubla2}
that the Lie algebroid structure of $\Om\to G$ may be prolonged to
$A\Om\to AG$.

That the anchor $\ts{a}\colon A\Om\to TA$ for the Lie algebroid of
$\Om\gpd A$ is a morphism of Lie algebroids over $a_G\colon AG\to TM$
follows as in the case of $A_2S$ above. The anchor ${\bf a}\colon A\Om\to TAG$
for the prolongation structure is $j_G^{-1}\circ A(\tila)$ where
$j_G\colon TAG\to ATG$ is the canonical isomorphism of
\cite[7.1]{MackenzieX:1994}. Since $\tila\colon\Om\to TG$ is a groupoid
morphism over $a\colon A\to TM$, and $j_G$ is an isomorphism of Lie
algebroids over $TM$, it follows that Condition~II is satisfied.
Condition~I is dealt with in the same way.

It was shown in \cite[\S3]{Mackenzie:SDGDPG} that $\Om^*\gpd K^*$, the
dual groupoid of $\Om$, together with the Poisson structures on $\Om^*$ and
$K^*$ dual to the Lie algebroid structures on $\Om$ and the core $K$, is a
Poisson groupoid. Thus Condition~III will follow from the next result.

\begin{thm}                               \label{thm:AOm}
The canonical isomorphism of double vector bundles $R = R^{gpd}\colon
A^*\Om^*\to A^*\Om = {\cal A}^{*V}$ is an isomorphism of Lie bialgebroids
$$
(A^*\Om^*, \Bar{A\Om^*})\to ({\cal A}^{*V}, ({\cal A}^{*V})^\dagger)
$$
where $(A\Om^*, A^*\Om^*)$ is the Lie bialgebroid of the Poisson
groupoid $\Om^*\gpd K^*$.
\end{thm}

We first recall the map $R$ from \cite[3.8]{Mackenzie:SDGDPG}. Associated
with $A\Om$ oriented as in Figure~\ref{fig:LAgpd}(b) there is the pairing
of the vertical and horizontal duals (\ref{eq:3duals}), which we write in
mnemonic form:
$$
\langle A^*\Om, A^\sol\Om\rangle_{K^*} =
\langle A^\sol\Om, A\Om\rangle_{AG} -
\langle A^*\Om, A\Om\rangle_{A}
$$
with the subscripts indicating the bases of the pairings. Using the
canonical isomorphism $I\colon A\Om^*\to A^\sol\Om$ induced by the
pairing $\llangle\ ,\ \rrangle\colon A\Om^*\pback{AG}\A\Om\to\R$,
we obtain, as in \cite[(18)]{Mackenzie:SDGDPG},
\begin{equation}                                   \label{eq:needI}
\ddagger A^*\Om, A\Om^*\ddagger =
\llangle A\Om^*, A\Om\rrangle -
\langle A^*\Om, A\Om\rangle_{A}.
\end{equation}
Now $R$ is defined by $\ddagger{\cal X}, R({\cal F})\ddagger =
\langle{\cal X}, {\cal F}\rangle$, where ${\cal X}\in A\Om^*,
{\cal F}\in A^*\Om^*$, and the pairing on the RHS is the standard one
over $K^*$. We finally have
\begin{equation}                                      \label{eq:3dualsAOm}
\langle{\cal X}, {\cal F}\rangle_{K^*} +
\langle R({\cal F}), \Xi\rangle_A  =
\llangle {\cal X}, \Xi\rrangle
\end{equation}
for compatible $\Xi\in A\Om$. Equivalently, the canonical isomorphism
$Z_V\colon (A^\sol\Om)^\dagger\to A^*\Om$ is given by
\begin{equation}                                  \label{eq:ZRI}
Z_V = R\circ I^\dagger.
\end{equation}

The first part of the following result was stated without proof in
\cite[\S3]{Mackenzie:SDGDPG}.

\begin{prop}                             \label{prop:RI}
{\em (i)} The map $R\colon A^*\Om^*\to A^*\Om$ is anti--Poisson from
the Poisson structure dual to the Lie algebroid of $\Om^*\gpd K^*$ to
the Poisson structure dual to the Lie algebroid of $\Om\gpd A$.

{\em (ii)} The map $I\colon A\Om^*\to A^\sol\Om$ is Poisson from the Poisson
structure induced on $A\Om^*$ {\em \cite{Weinstein:1988}} by the Poisson
groupoid structure on $\Om^*\gpd K^*$, to the Poisson structure dual to
the prolonged Lie algebroid structure on $A\Om\to AG$.
\end{prop}

\pf
It suffices \cite{Weinstein:1988} to prove that the graph of $R$ is
coisotropic in $A^*\Om^*\times A^*\Om$. Let $S = \Om^*\pback{G}\Om$
and write $F\colon S\to\R$ for the pairing. Then $F$ is a groupoid
morphism, where $S$ is the pullback groupoid over $K^*\pback{M}A$,
and so, as in \cite[3.7]{Mackenzie:SDGDPG}, we can apply the Lie
functor and get $\llangle\ ,\ \rrangle = A(F)\colon AS\to\R$. This is
linear and so defines a section $\nu$ of the dual of $AS$, which is closed
since $A(F)$ is a morphism. So by \cite[4.6]{MackenzieX:1998}, the image
of $\nu$ is coisotropic.

It remains to show that the image of $\nu$ coincides with the graph
of $R$. The image of $\nu$ consists of those $({\cal F}, {\cal X})\in
A^*_\kappa\Om\times A^*_Y\Om$ such that
$$
\langle({\cal F}, \Phi\rangle, ({\cal X}, \Xi)\rangle =
A(F)({\cal X}, Xi)
$$
for all $({\cal X}, \Xi)\in AS$ compatible with $({\cal F}, \Phi)$.
As in \cite[5.5]{MackenzieX:1994}, this equation expands to
(\ref{eq:3dualsAOm}).

We leave the proof of (ii) to the reader.
\boom

\noindent
{\sc Proof of Theorem \ref{thm:AOm}:} We must first show that $R$ is an
isomorphism of Lie algebroids $A^*\Om^*\to {\cal A}^{*V}$. Now the Lie
algebroid structure on ${\cal A}^{*V}$ is induced from
$({\cal A}^{*H})^\dagger$ via $Z_V$. So what we have to show is that
$Z_V^{-1}\circ R\colon A^*\Om^*\to (A^\sol\Om)^\dagger$ is an isomorphism
of Lie algebroids, and this is equivalent to the dual over $K^*$ being
Poisson. This dual is, using (\ref{eq:ZRI}), $I^{-1}\colon A^\sol\Om\to
A\Om^*$, and so the result follows from \ref{prop:RI}(ii) above.

Secondly we must show that
$R^\dagger\colon ({\cal A}^{*V})^\dagger\to\Bar{A\Om^*}$ is an isomorphism
of Lie algebroids over $K^*$. (Note that the minus sign is in the
bundle over $K^*$.) This is equivalent to showing that the dual
$R\colon\Bar{A^*\Om^*}\to A^*\Om$ is Poisson, and this is
\ref{prop:RI}(i) above.
\boom

In summary, we have proved:

\begin{thm}
The double Lie algebroid $(A\Om;A,AG;M)$ of an \LAgpd\
$(\Om;A,G;M)$, as constructed in {\em \cite[\S1]{Mackenzie:Doubla2}}, is a
double Lie algebroid as defined in {\em \ref{df:doubla}}.

In particular, the double Lie algebroids $(A_2S;AH,AV;M)$ and
$(A^2S;AH,AV;M)$ of a double Lie groupoid $(S;H,V;M)$, as constructed
in {\em \cite[\S2]{Mackenzie:Doubla2}}, are double Lie algebroids as
defined in {\em \ref{df:doubla}}.
\end{thm}

\begin{ex}\rm                                          \label{ex:TA}
Let $A$ be any Lie algebroid on $M$. Then
$\Om = A\times A$ has an \LAgpd\ structure over $M\times M$ and $A$,
and the associated double Lie algebroid constructed in
\cite[\S1]{Mackenzie:Doubla2} is $(TA;A,TM;M)$.

The associated duals are ${\cal A}^{*V} = T^*A$ and
${\cal A}^{*H} = T^\sol A$. Using $R$ and $I$ as in \cite{MackenzieX:1994},
these can be identified with $T^*A^*$ and $T(A^*)$, as bundles over $A^*$.
The Lie algebroid structure on $T^*A^*$ is the cotangent of the dual Poisson
structure on $A^*$. The Lie algebroid structure on $T(A^*)$ is the standard
tangent bundle structure. This is the standard Lie bialgebroid $(T^*P, TP)$
for $P = A^*$.
\end{ex}

\begin{ex}\rm
Taking $A = TM$ in the previous example, we see that $T^2M$ is a double
Lie algebroid with associated bialgebroid $(T^*T^*M,TT^*M)$. This is a
Lie bialgebroid over $T^*M$, the induced Poisson structure being the
standard symplectic structure.
\end{ex}

The double Lie algebroids considered in the next two sections do not
necessarily have an underlying \LAgpd.

\section{THE DOUBLE LIE ALGEBROID OF A LIE \break BIALGEBROID}
\label{sect:dlalba}

Here we use the following criterion for a Lie bialgebroid.

\begin{thm}                                        \label{thm:6.2}
{\bf \cite[6.2]{MackenzieX:1994}}
Let $A$ be a Lie algebroid on $M$ such that its dual vector bundle $A^{*}$
also has a Lie algebroid structure. Denote their anchors by $a, a_*$.
Then $(A, A^{*})$ is a Lie bialgebroid if and only if
\begin{equation}                         \label{main}
T^*(A^*)\buildrel{R}\over\longrightarrow T^*(A)
\buildrel\pi^\#_A\over\longrightarrow TA
\end{equation}
is a Lie algebroid morphism over $a_*$, where the domain
$T^{*}(A^{*})\to A^{*}$ is the cotangent Lie algebroid induced by the
Poisson structure on $A^{*}$, and the target $TA\to TM$ is the tangent
prolongation of $A$.
\end{thm}

Consider a Lie algebroid $A$ on $M$ together with a Lie algebroid
structure on the dual, not {\em a priori} related to that on $A$. The
structure on $A^*$ induces a Poisson structure on $A$, and this gives
rise to a cotangent Lie algebroid $T^*A\to A$. Equally, the Lie algebroid
structure on $A$ induces a Poisson structure on $A^*$ and this gives rise
to a cotangent Lie algebroid $T^*A^*\to A^*$. We transfer this latter
structure to $T^*A\to A^*$ via $R$.

There are now four Lie algebroid structures on the four sides of
${\cal A} = T^*A$ as in Figure~\ref{fig:TA}(b).

\begin{thm}                                       \label{thm:Manin}
Let $A$ be a Lie algebroid on $M$ such that its dual vector bundle $A^{*}$
also has a Lie algebroid structure. Then $(A, A^*)$ is a Lie bialgebroid
if and only if ${\cal A} = T^*A$, with the structures just described, is a
double Lie algebroid.
\end{thm}

\pf
Assume that $(A, A^*)$ is a Lie bialgebroid. The vertical structure
on ${\cal A}$ is the cotangent Lie algebroid structure for the Poisson
structure on $A$. The anchor of this is a morphism of double vector
bundles $\pi^\#_A\colon T^*A\to TA$ over $a_*\colon A^*\to TM$ and $\id_A$,
inducing $-a^*_*\colon T^*M\to A$ on the cores. Now the horizontal structure
has the cotangent Lie algebroid structure for the Poisson structure on
$A^*$, transported via $R = R_A\colon T^*A^*\to T^*A$. So the condition that
$\pi^\#_A$ is a morphism of Lie algebroids over $a_*$ with respect to the
horizontal structure is precisely \ref{thm:6.2}.

On the other hand, the anchor for the horizontal structure is
$$
\pi^\#_{A^*}\circ R^{-1}\colon T^*A\to T(A^*),
$$
and this is a morphism of double vector bundles over $a\colon A\to TM$
and $\id_{A^*}$, inducing $+a^*$ on the cores. Since $R^{-1} = R_{A^*}$,
the condition that this anchor be a morphism with respect to the
vertical structure is precisely the dual form of \ref{thm:6.2}, to which
\ref{thm:6.2} is equivalent by \cite[3.10]{MackenzieX:1994} or
\cite{Kosmann-Schwarzbach:1995}.

The vertical dual of ${\cal A}$ is the tangent double vector bundle as in
Figure~\ref{fig:TA}(a). Being the dual of a Lie algebroid, the vertical
structure of this has a Poisson structure; this is the tangent lift
of the Poisson structure on $A$ \cite[5.6]{MackenzieX:1994}. The
corresponding Poisson tensor is
$$
\pi^\#_{TA} = J_A\circ T(\pi^\#_A)\circ\theta_A^{-1}
$$
where $J_A\colon T^2A\to T^2A$ is the canonical involution for the
manifold $A$ and
$\theta_A\colon T(T^*A)\to T^*(TA)$ is the canonical map $\alpha$ of
\cite{Tulczyjew}, denoted $J'$ in \cite[5.4]{MackenzieX:1994}.

We must check that this Poisson structure coincides with that induced from
the dual of ${\cal A}^{*H}$. Note that we have $K^* = TM$ here and to avoid
confusion we drop the ${}^\dagger$ notation in this case and denote all
duals over $TM$ by ${}^\sol$. Duals over $A^*$ will be denoted
${}^{*A^*}$.

Consider $I\circ R^{*A^*}\colon (T^*A)^{*A^*}\to T^\sol A$. This preserves
the core $A^*$ and the side $A^*$ but reverses the side $TM$. Define
$W = -I\circ R^{*A^*}$ where the minus is for the bundle over $A^*$.
Then $W^\sol\colon TA\to ({\cal A}^{*H})^\sol$ and the reader can
check that $Z_V^{-1} = \dminus W^\sol$ where the heavy minus is over $TM$.
See alternatively \cite[3.3]{Mackenzie:SDGDPG}.

To prove that $Z_V\colon ({\cal A}^{*V})^\sol \to TA$ is an isomorphism
of Lie algebroids over $TM$ we must show that $W$ is an anti--Poisson
map. This may be done directly or by observing that $-W$ is, in terms
of the double Lie algebroid ${\cal A}' = TA$ of \ref{ex:TA}, the map
$(Z'_V)^\dagger = Z'_H$.

So we have ${\cal A}^{*V} = TA$ and $({\cal A}^{*V})^\sol = T^\sol A$ and
Condition~III follows from the next result. We use $I$ to replace $T^\sol A$
by $TA^*$.

\begin{lem}
Given that $(A,A^*)$ is a Lie bialgebroid on $M$, the tangent prolongation
structures make $(TA, TA^*)$ a Lie bialgebroid on $TM$ with respect to the
tangent pairing.
\end{lem}

\pf
We use the bialgebroid criterion of \ref{thm:6.2}. We must
prove that
\begin{equation}                                \label{eq:comp}
T^*(TA^*)\buildrel\Tilde{R}\over\longrightarrow T^*(TA)
\buildrel\pi^\#_{TA}\over\longrightarrow T^2A
\end{equation}
is a morphism of Lie algebroids over the anchor of $TA^*$ which, by
\cite[5.1]{MackenzieX:1994}, is $J_M\circ T(a_*)\colon TA^*\to T^2M$.
Here $\Tilde{R}$ is the canonical map $R$ for $TA\to TM$, transported
using $I\colon TA^*\to T^\sol A$. The domain of (\ref{eq:comp}) is the
cotangent Lie algebroid for the Poisson structure on $TA^*$, which Poisson
structure---again by \cite[5.6]{MackenzieX:1994}---is both the tangent lift
of the Poisson structure on $A^*$ and the dual (via $I$) of the prolongation
Lie algebroid structure on $TA$. The target of (\ref{eq:comp}) is the
iterated tangent prolongation of the Lie algebroid structure of $A$.

Now $\Tilde{R} = \theta_A\circ T(R_A)\circ\theta_{A^*}^{-1}$ and so
$$
\pi^\#_{TA}\circ \Tilde{R} =
J_A\circ T(\pi^\#_A\circ R_A)\circ\theta_{A^*}^{-1}.
$$
We know that $\pi^\#_A\circ R_A\colon T^*A^*\to TA$ is a morphism of
Lie algebroids over $a_*$, so $T(\pi^\#_A\circ R_A)$ is a morphism of
the prolongation structures over $T(a_*)$. We need two further
observations.

Firstly, for any Poisson manifold, $\theta_P\colon T(T^*P)\to T^*(TP)$
is an isomorphism of Lie algebroids over $TP$ from the tangent prolongation
of the cotangent Lie algebroid structure on $T^*P$ to the cotangent Lie
algebroid of the tangent Poisson structure \cite[2.13]{Mackenzie:Doubla2}.
We apply this to $P = A^*$.

Secondly, $J_A\colon T^2A\to T^2A$ is a Lie algebroid automorphism over
$J_M$ of the iterated prolongation of the given Lie algebroid
structure on $A$.

Putting these facts together, we have that $\pi^\#_{TA}\circ\Tilde{R}$ is a
Lie algebroid morphism.
\boom

Now conversely suppose that $A$ is a Lie algebroid on $M$ and that $A^*$
has a Lie algebroid structure, not {\em a priori} related to the structure
on $A$. Consider ${\cal A} = T^*A$ with the two cotangent Lie algebroid
structures arising from the Poisson structures on $A^*$ and $A$, and
suppose that these structures make ${\cal A}$ a double Lie algebroid.

Then in particular the anchor of the horizontal structure must be a Lie
algebroid morphism with respect to the other structures, as in Condition~II,
and this is
$$
T^*A\buildrel{R_{A^*}}\over\longrightarrow T^*A^*
\buildrel\pi^\#_{A^*}\over\longrightarrow TA^*
$$
That this be a Lie algebroid morphism over $a\colon A\to TM$ is
precisely the dual form of \ref{thm:6.2}.

This completes the proof of Theorem \ref{thm:Manin}.
\boom

Recall the Manin triple characterization of a Lie bialgebra, as given
in \cite{LuW:1990}: Given a Lie bialgebra $({\goth g}, {\goth g}^*)$
the vector space direct sum ${\goth d} = {\goth g}\oplus{\goth g}^*$
has a Lie algebra bracket defined in terms of the two coadjoint
representations. This bracket is invariant under the pairing
$\langle X + \phi, Y + \psi\rangle = \psi(X) + \phi(Y)$ and both ${\goth g}$
and ${\goth g}^*$ are coisotropic subalgebras. Conversely, if a Lie
algebra ${\goth d}$ is a vector space direct sum ${\goth g}\oplus{\goth h}$,
both of which are coisotropic with respect to an invariant pairing
of ${\goth d}$ with itself, then ${\goth h}\isom{\goth g}^*$ and
$({\goth g}, {\goth h})$ is a Lie bialgebra, with ${\goth d}$ as the
double.

Two aspects of this result concern us here. Firstly, it characterizes the
notion of Lie bialgebra in terms of a single Lie algebra structure on
${\goth d}$, the conditions being expressed in terms of the simple notion
of pairing. Secondly, the roles of the two Lie algebras ${\goth g}$ and
${\goth g}^*$ are completely symmetric; it is an immediate consequence of
the Manin triple result that $({\goth g}, {\goth g}^*)$ is a Lie bialgebra
if and only if $({\goth g}^*, {\goth g})$ is so.

In considering a corresponding characterization for Lie algebroids, the
most important difference to note is that the structure on the double is
no longer over the same base as the given Lie algebroids. This is to be
expected in view of the results of \cite[\S2]{Mackenzie:1992} for the
double groupoid case. There it is proved that if a double groupoid
$(S;H,V;M)$ has trivial core (that is, the only elements of $S$ to have
two touching sides which are identity elements, are those which are
identities for both structures), then there is a third groupoid structure
on $S$, over base $M$, called in \cite[p.200]{Mackenzie:1992} the
{\em diagonal structure} and denoted $S_D$. With respect to this structure
on $S$, the identity maps from $H$ and $V$ into $S$ are morphisms over
$M$, and $S$ as a manifold is $H\pback{M}V$. This diagonal structure is, in
the case where $H$ and $V$ are dual Poisson groups, precisely the structure
of the double group. The existence of the diagonal structure, in the general
formulation given in \cite[\S2]{Mackenzie:1992}, depends crucially on the
fact that $S$ has trivial core; that is, that $S$ is vacant. Since the core
of the double vector bundle $T^*A$, for $A$ a vector bundle on $M$, is
$T^*M$, we expect $T^*A$ to possess a Lie algebroid structure over $M$ only
when $M$ is a point.

The role played in the bialgebra case by the Lie algebra structure of
${\goth d} = {\goth g}\oplus{\goth g}^*$ is taken, for Lie bialgebroids,
by the two structures on $T^*A$ (whose bases are $A$ and $A^*$). In place
of a characterization in terms of a single Lie algebroid structure on
$T^*A$ with base $M$, Theorem \ref{thm:Manin} gives a characterization in
terms of the two Lie algebroid structures on $T^*A$. The role of the pairing
in the bialgebra case is taken in \ref{thm:Manin} by the isomorphism $R$.

The analysis of Lie bialgebras and Poisson Lie groups is usually given
in terms of the coadjoint representations and the dressing transformation
actions. This was extended by \cite{LuW:1990}, \cite{Majid:1990} and
\cite{Kosmann-SchwarzbachM:1988} to the more general situation of matched
pairs of Lie groups and Lie algebras, and in \cite{Mackenzie:1992} to matched
pairs of groupoids. In the next section we consider the corresponding results
for Lie algebroids.

\section{MATCHED PAIRS AND VACANT DOUBLE LIE \break ALGEBROIDS}
\label{sect:mpvdla}

The notion of matched pair of Lie algebras was introduced by
Kosmann--Schwarzbach and Magri \cite{Kosmann-SchwarzbachM:1988}, who
called them {\em extensions bicrois\'ees}, by Lu and Weinstein
\cite{LuW:1990}, who called them {\em double Lie algebras}, and by
Majid \cite{Majid:1990}, who introduced the term {\em matched pair}.
(In fact, forms of the concept had been found much earlier; see
\cite{Weinstein:1990} for references.) A matched pair of Lie algebras may
be regarded as a triple $({\goth d}, {\goth g}, {\goth h})$ where the Lie
algebra ${\goth d}$ is the vector space direct sum of its subalgebras
${\goth g}$ and ${\goth h}$; a matched pair is thus the notion of Manin
triple with the duality aspect removed. A matched pair can be described
in terms of a pair of representations, of ${\goth g}$ on ${\goth h}$
and of ${\goth h}$ on ${\goth g}$, subject to twisted derivation
conditions. See the references above or \ref{df:mokri}, \ref{prop:mokri}
below.

The corresponding concept of matched pair of Lie groups \cite{LuW:1990},
\cite{Majid:1990} was extended to groupoids
in \cite{Mackenzie:1992}. In \cite{Mokri:1997}, Mokri differentiated the
twisted automorphism equations of \cite{Mackenzie:1992} to obtain conditions
on a pair of Lie algebroid representations, of $A$ on $B$ and of $B$ on $A$,
which ensure that the direct sum vector bundle $A\oplus B$ has a Lie
algebroid structure with $A$ and $B$ as subalgebroids. We quote the
following.

\begin{df}  {\bf \cite[4.2]{Mokri:1997}}   \label{df:mokri}
Let $A$ and $B$ be Lie algebroids on base $M$, with anchors $a$ and $b$,
and let $\rho\colon A\to\CDO(B)$ and $\sigma\colon B\to\CDO(A)$ be
representations of $A$ on the vector bundle $B$ and of $B$ on the vector
bundle $A$. Then $A$ and $B$ together with $\rho$ and $\sigma$ form a
{\em matched pair} if the following equations hold for all $X, X_1, X_2
\in \Ga A,\ Y, Y_1, Y_2\in\Ga B$:
\begin{eqnarray*}
\rho_X([Y_1, Y_2]) & = & [\rho_X(Y_1), Y_2] + [Y_1, \rho_X(Y_2)]
   +\rho_{\sigma_{Y_2}(X)}(Y_1) - \rho_{\sigma_{Y_1}(X)}(Y_2),\\
\sigma_Y([X_1, X_2]) & = & [\sigma_Y(X_1), X_2] + [X_1, \sigma_Y(X_2)]
   +\sigma_{\rho_{X_2}(Y)}(X_1) - \sigma_{\rho_{X_1}(Y)}(X_2),\\
a(\sigma_Y(X)) - b(\rho_X(Y)) & = & [b(Y), a(X)].
\end{eqnarray*}
\end{df}

Here $\CDO(E)$, for any vector bundle $E$, is the vector bundle whose
sections are the first or zeroth order differential operators
$D\colon\Ga E\to\Ga E$ for which there is a vector field $X$ on $M$ with
$D(f\mu) = fD(\mu) + X(f)\mu$ for all $f\in C(M),\ \mu\in\Ga E$. With
anchor $D\mapsto X$ and the usual bracket, $\CDO(E)$ is a Lie algebroid
(see \cite[III\S2]{Mackenzie:LGLADG}).

\begin{prop} {\bf \cite[4.3]{Mokri:1997}}        \label{prop:mokri}
Given a matched pair, there is a Lie algebroid structure on the direct sum
vector bundle $A\oplus B$, with anchor $c(X\oplus Y) = a(X) + b(Y)$ and
bracket
\begin{equation}
[X_1\oplus Y_1, X_2\oplus Y_2] =
   \{[X_1, X_2] + \sigma_{Y_1}(X_2) - \sigma_{Y_2}(X_1)\}\oplus
   \{[Y_1, Y_2] + \rho_{X_1}(Y_2) - \rho_{X_2}(Y_1)\}.
\end{equation}
Conversely, if $A\oplus B$ has a Lie algebroid structure making
$A\oplus 0$ and $0\oplus B$ Lie subalgebroids, then $\rho$ and $\sigma$
defined by $[X\oplus 0, 0\oplus Y] = -\sigma_Y(X) \oplus \rho_X(Y)$
form a matched pair.
\end{prop}

We now show that matched pairs correspond precisely to double Lie algebroids
with zero core. The following definition is a natural sequel to
\cite[2.11, 4.10]{Mackenzie:1992}.

\begin{df}
A double Lie algebroid $({\cal A}; A^H, A^V; M)$ is {\em vacant} if the
combination of the two projections,
$(\tilq_V, \tilq_H)\colon {\cal A}\to A^H\pback{M} A^V$ is a diffeomorphism.
\end{df}

Consider a vacant double Lie algebroid, which we will write here as
$({\cal A}; A, B; M)$. Note that ${\cal A}\to A^H$ and ${\cal A}\to A^V$
are the pullback bundles $q_A^!B$ and $q_B^!A$. The two duals are
${\cal A}^{*H} = A^*\oplus B$ and ${\cal A}^{*V} = A\oplus B^*$, as vector
bundles over $M$, and the duality is (see \cite[3.4]{Mackenzie:SDGDPG})
\begin{equation}                                         \label{eq:pairing}
\langle X + \psi, \phi + Y\rangle =
                           \langle \phi, X \rangle - \langle\psi, Y\rangle.
\end{equation}

The horizontal bundle projection $\tilq_A\colon {\cal A}\to B$ is a
morphism of Lie algebroids over $q_A\colon A\to M$ and since it is a
fibrewise surjection, it defines an action of $B$ on $q_A$ as in
\cite[\S2]{HigginsM:1990a}. Namely, each section $Y$ of $B$ induces the
pullback section $1\otimes Y$ of $q_A^!B$ and this induces a vector field
$\eta(Y) = \tilb(1\otimes Y)$ on $A$, where $\tilb\colon{\cal A}\to TA$ is
the anchor of the vertical structure. By Conditions I and II, $\eta(Y)$ is
linear over the vector field $b(Y)$ on $M$, in the sense of
\cite[\S1]{MackenzieX:1998}; that is, $\eta(Y)$ is a vector bundle morphism
$A\to TA$ over $b(Y)\colon M\to TM$. It follows that $\eta(Y)$ defines
covariant differential operators $\sigma^{(*)}_Y$ on $A^*$ and $\sigma_Y$
on $A$ by
\begin{equation}                          \label{eq:sigma}
\eta(Y)(\ell_\phi) = \ell_{\sigma^{(*)}_Y(\phi)},\quad
\langle\phi, \sigma_Y(X)\rangle = b(Y)\langle\phi, X\rangle
- \langle\sigma^{(*)}_Y(\phi), X\rangle
\end{equation}
where $\phi\in\Ga A^*,\ X\in\Ga A,$ and $\ell_\phi$ denotes the function
$A\to\R,\ X\mapsto\langle\phi(q_AX), X\rangle$; see
\cite[\S2]{MackenzieX:1998}. Since $\tilq_A$ is a Lie algebroid
morphism, it follows that $\sigma$ is a representation of $B$ on the
vector bundle $A$.

Dually, $\tilq_B$ is a morphism of Lie algebroids over $q_B$ and for
each $X\in\Ga A$ we obtain a linear vector field $\xi(X)\in {\cal X}(B)$
over $a(X)$. We likewise define covariant differential operators
$\rho^{(*)}_X$ on $B^*$ and $\rho_X$ on $B$ by
\begin{equation}                             \label{eq:rho}
\xi(X)(\ell_\psi) = \ell_{\rho^{(*)}_X(\psi)},\quad
\langle\psi, \rho_X(Y)\rangle = a(X)\langle\psi, Y\rangle
- \langle\rho^{(*)}_X(\psi), Y\rangle.
\end{equation}
Again, $\rho^{(*)}_X$ and $\rho_X$ are representations of $A$.

In fact (see \cite[\S2]{HigginsM:1990a}) the two Lie algebroid structures
on ${\cal A}$ are action Lie algebroids determined by the actions
$Y\mapsto\eta(Y)$ and $X\mapsto\xi(X)$. It follows that the dual Poisson
structures are semi--direct in a general sense, but we prefer to proceed
on an ad hoc basis.

For a general vector bundle, the functions on the dual are generated by
the linear functions and the pullbacks from the base manifold. In the case
of a pullback bundle such as $q_B^!A$, one can refine this description a
little further. Namely, if $\pi_B\colon q_B^!A^*\to A^*$ is
$(\phi, Y)\mapsto \phi$, and $\tilq_{A^*}\colon q_B^!A^*\to B$ is the bundle
projection, then the functions on $q_B^!A^*$ are generated by all
$$
\ell_X\circ\pi_B,\qquad
\ell_\psi\circ\tilq_{A^*}\qquad\mbox{and}\qquad
f\circ q_B\circ\tilq_{A^*}
$$
where $X\in\Ga A,\ \psi\in\Ga B^*$ and $f\in C(M)$. Now the Poisson
structure on $q_B^!A^*$ is characterized by
\begin{eqnarray}                              \label{eq:poissA*}
\{\ell_{X_1}\circ\pi_B,\ \ell_{X_2}\circ\pi_B\}
   = \ell_{[X_1,X_2]}\circ\pi_B, & \qquad &
\{\ell_{X}\circ\pi_B,\ \ell_\psi\circ\tilq_{A^*}\}
   = \ell_{\rho^{(*)}_X(\psi)}\circ \tilq_{A^*},\nonumber \\
\{\ell_{X}\circ\pi_B,\ f\circ q_B\circ\tilq_{A^*}\}
   = a(X)(f)\circ q_B\circ \tilq_{A^*}, & \qquad &
\{F_1\circ\tilq_{A^*},\ F_2\circ\tilq_{A^*}\}
   = 0,
\end{eqnarray}
where $F_1, F_2$ are any smooth functions on $B$. Similarly,
\begin{eqnarray}                            \label{eq:poissB*}
\{\ell_{Y_1}\circ\pi_A,\ \ell_{Y_2}\circ\pi_A\}
   = \ell_{[Y_1,Y_2]}\circ\pi_A, & \qquad &
\{\ell_{Y}\circ\pi_A,\ \ell_\phi\circ\tilq_{B^*}\}
   = \ell_{\sigma^{(*)}_Y(\phi)}\circ \tilq_{B^*},\nonumber \\
\{\ell_{Y}\circ\pi_A,\ f\circ q_A\circ\tilq_{B^*}\}
   = b(Y)(f)\circ q_A\circ \tilq_{B^*}, & \qquad &
\{G_1\circ\tilq_{B^*},\ G_2\circ\tilq_{B^*}\}
   = 0,
\end{eqnarray}

Now these Poisson structures induce Lie algebroid structures on the
direct sum bundles $A\oplus B^*$ and $A^*\oplus B$ over $M$. Consider
first a section $X\oplus\psi$ of $\Ga(A\oplus B^*)$. Via the pairing
(\ref{eq:pairing}), this induces a linear function on $q_B^!A^*$, namely
$$
\ell^\dagger_{X\oplus\psi} = \ell_X\circ\pi_B - \ell_\psi\circ\tilq_{A^*}
$$
where $\ell^\dagger$ refers to the pairing (\ref{eq:pairing}). By following
through the equations (\ref{eq:poissA*}) and (\ref{eq:poissB*}) one
obtains the following.

\begin{lem}
The Lie algebroid structure on $A\oplus B^*$ induced as above has anchor
$e(X\oplus\psi) = a(X)$ and bracket
\begin{equation}                            \label{eq:sdprho}
[X_1\oplus\psi_1, X_2\oplus\psi_2] =
   [X_1, X_2]\oplus\{\rho^{(*)}_{X_1}(\psi_2) - \rho^{(*)}_{X_2}(\psi_1)\}.
\end{equation}

The Lie algebroid structure on $A^*\oplus B$ induced as above has anchor
$e_*(\phi\oplus Y) = -b(Y)$ and bracket
\begin{equation}                              \label{eq:sdpsigma}
[\phi_1\oplus Y_1, \phi_2\oplus Y_2] =
   \{\sigma^{(*)}_{Y_2}(\phi_1) - \sigma^{(*)}_{Y_1}(\phi_2)\}\oplus[Y_2, Y_1].
\end{equation}
\end{lem}

Thus $A\oplus B^*$ is the semi--direct product (over the base $M$, in the
sense of \cite{Mackenzie:LGLADG}) of $A$ with the vector bundle $B^*$
with respect to $\rho^{(*)}$. However $A^*\oplus B$ is the opposite Lie
algebroid to the semi--direct product of $B$ with $A^*$.

We can now apply Condition III to $A\oplus B^*$ and $A^*\oplus B$. For
brevity write $E = A\oplus B^*$. Recall \cite{MackenzieX:1994},
\cite{Kosmann-Schwarzbach:1995} that $(E, E^*)$ is a Lie bialgebroid if
and only if
\begin{equation}                         \label{eq:full}
d^E[\phi_1 \oplus Y_1, \phi_2\oplus Y_2] =
   [d^E(\phi_1\oplus Y_1), \phi_2\oplus Y_2] +
   [\phi_1\oplus Y_1, d^E(\phi_2\oplus Y_2)]
\end{equation}
for all $\phi_1\oplus Y_1, \phi_2\oplus Y_2\in \Ga E^*$. It follows
\cite[3.4]{MackenzieX:1994} that for any $f\in C(M),\ X\oplus\psi\in\Ga E$,
\begin{equation}                     \label{eq:lied}
L_{d^Ef}(X\oplus\psi) + [d^{E^*}f, X\oplus\psi] = 0.
\end{equation}
It is easy to check that
$$
d^Ef = d^Af\oplus 0,\qquad d^{E^*}f = 0\oplus d^Bf;
$$
note that these imply that the Poisson structure induced on $M$ by the Lie
bialgebroid $(E, E^*)$ is zero. Now the Lie derivative in (\ref{eq:lied})
is a standard Lie derivative for $E^*$ and so
\begin{eqnarray}                        \label{eq:shift1}
\langle L_{d^Ef}(X\oplus\psi), \phi\oplus Y\rangle
   & = & e_*(d^Ef)\langle X\oplus\psi, \phi\oplus Y\rangle
         -\langle X\oplus\psi, [d^Ef, \phi\oplus Y]\rangle \nonumber\\
   & = & 0 - \langle X\oplus\psi, \sigma^{(*)}_Y(d^Af)\oplus 0\rangle \nonumber\\
   & = & \langle d^Af, \sigma_Y(X)\rangle - b(Y)\langle d^Af, X\rangle \nonumber\\
   & = & a(\sigma_Y(X))(f) - b(Y)a(X)(f)
\end{eqnarray}
where we used (\ref{eq:sigma}) and (\ref{eq:sdpsigma}). Expanding out the
bracket term in (\ref{eq:lied}) in a similar way, we obtain the third
equation in \ref{df:mokri}.

Now consider the bialgebroid equation (\ref{eq:full}). We calculate this
in the case $\phi_1 = \phi_2 = 0$, with arguments $0\oplus \psi_1,\
X_2\oplus 0.$ With these values we refer to (\ref{eq:full}) as equation
(\ref{eq:full}$)_0$. First we need the following lemma, which is a
straightforward calculation.

\begin{lem}
$$
d^E(\phi\oplus Y)(X_1\oplus\psi_1, X_2\oplus\psi_2) =
   (d^A\phi)(X_1, X_2) + \langle\psi_1, \rho_{X_2}(Y)\rangle
            - \langle\psi_2, \rho_{X_1}(Y)\rangle.
$$
\end{lem}

The LHS of (\ref{eq:full}$)_0$ is easily seen to be
$\langle\psi_1, \rho_{X_2}[Y_2, Y_1]\rangle.$ On the RHS, consider
the second term first. Regarding the bracket as a Lie derivative, we have
\begin{eqnarray}                    \label{eq:shift2}
\lefteqn{\langle L_{0\oplus Y_1}(d^E(0\oplus Y_2)),
   (0\oplus\psi_1)\wedge(X_2\oplus 0)\rangle}\\
& = & L_{0\oplus Y_1}\langle d^E(0\oplus Y_2), (0\oplus\psi_1)\wedge(X_2\oplus 0)\rangle
- \langle d^E(0\oplus Y_2), L_{0\oplus Y_1}((0\oplus\psi_1)\wedge(X_2\oplus 0))\rangle
\nonumber
\end{eqnarray}
Since $e_*(0\oplus Y_1) = -b(Y_1)$, the first term is
$-b(Y_1)\langle\psi_1, \rho_{X_2}(Y_2)\rangle.$ For the second term we
need the following lemma.

\begin{lem}                              \label{lem:lied}
For any $\phi\oplus Y\in \Ga E^*,\ X\oplus\psi\in\Ga E$, we have
$L_{\phi\oplus Y}(X\oplus\psi) = -\sigma_Y(X)\oplus \Bar{\psi}$ where
for any $Y'\in B$,
$$
\langle\Bar{\psi}, Y'\rangle = -b(Y)\langle\psi, Y'\rangle
   +\langle\sigma^{(*)}_{Y'}(\phi), X\rangle + \langle\psi, [Y,Y']\rangle.
$$
\end{lem}

\pf
This is a Lie derivative for $E^*$ and applying the same device as in
(\ref{eq:shift1}) we have, for any $\phi\oplus Y'\in\Ga E^*$,
$$
\langle\phi'\oplus Y', L_{\phi\oplus Y}(X\oplus\psi)\rangle
= -\langle\phi', \sigma_Y(X)\rangle + b(Y)\langle\psi, Y'\rangle
-\langle\sigma^{(*)}_{Y'}(\phi), X\rangle + \langle\psi, [Y',Y]\rangle.
$$
Setting $Y' = 0$ and $\phi' = 0$ in turn gives the result.
\boom

Now expand out the Lie derivative of the wedge product in (\ref{eq:shift2})
and apply Lemma \ref{lem:lied}. One obtains for the second term on the
RHS of (\ref{eq:full}$)_0$
$$
-\langle\psi_1, [Y_1, \rho_{X_2}(Y_2)]\rangle
+ \langle\psi_1, \rho_{\sigma_{Y_1}(X_2)}(Y_2)\rangle.
$$
The first term on the RHS of (\ref{eq:full}$)_0$ is easily obtained from
this, and combining with the LHS we have the first equation in
\ref{df:mokri}. The second equation is obtained in a similar way from the
dual form of (\ref{eq:full}).

This completes the proof of the first part of the following result. The
second part is proved essentially by reversing the steps.

\begin{thm}                                     \label{thm:mp}
Let $({\cal A};A,B;M)$ be a vacant double Lie algebroid. Then the two
Lie algebroid structures on ${\cal A}$ are action Lie algebroids
corresponding to actions which define representations $\rho$, of $A$
on $B$, and $\sigma$, of $B$ on $A$, with respect to which $A$ and $B$
form a matched pair.

Conversely, let $A$ and $B$ be a matched pair of Lie algebroids with
respect to $\rho$ and $\sigma$. Then the action of $A$ on $q_B$ induced
by $\rho$ and the action of $B$ on $q_A$ induced by $\sigma$ define Lie
algebroid structures on ${\cal A} = A\pback{M}B$ with respect to which
$({\cal A};A,B;M)$ is a vacant double Lie algebroid.
\end{thm}

In the case of a Lie bialgebra $({\goth g}, {\goth g}^*)$, the Lie
bialgebroid associated to the vacant double Lie algebroid is thus
$({\goth g}\oplus{\goth g}_0, {\goth g}^*\oplus{\goth g}^*_0)$ where the
subscripts denote the abelianizations. This is of course consistent with
\ref{thm:Manin} in the bialgebra case---which is both a bialgebroid and a
matched pair.

One other example which should be mentioned briefly is that of affinoids.
An {\em affinoid} \cite{Weinstein:1990} may be regarded as a vacant double
Lie groupoid in which both side groupoids $H$ and $V$ are the graphs of
simple foliations defined by surjective submersions $\pi_1\colon M\to Q_1$
and $\pi_2\colon M\to Q_2$. The corresponding double Lie algebroid was
calculated in \cite{Mackenzie:Doubla2} to be a pair of conjugate flat
partial connections adapted to the two foliations. The bialgebroid in
this case is $(T^{\pi_1}M\oplus(T^{\pi_2}M)^*,
(T^{\pi_1}M)^*\oplus T^{\pi_2}M)$ with semi--direct structures defined
by the connections.

Theorem \ref{thm:mp} provides a diagrammatic characterization of matched
pairs of Lie algebroids, directly comparable to the diagrammatic description
of matched pairs of group(oid)s given in \cite[\S2]{Mackenzie:1992}. In the
groupoid case, the twisted multiplicativity equations are fairly unintuitive,
and we believe that the derivation of them directly from the vacant double
groupoid axioms has been a significant clarification. In the Lie algebroid
case the equations in \ref{df:mokri} are again not simple and, unlike the
groupoid case, are defined in terms of sections rather than elements.
Nonetheless the characterization given by \ref{thm:mp} is purely diagrammatic:
recall that the characterization \ref{thm:6.2} of a Lie bialgebroid is
formulated entirely in terms of the Poisson tensor and the canonical
isomorphism $R$. Thus we have a definition of matched pair which can be
formulated more generally in a category possessing pullbacks and suitable
additive structure.

We will show elsewhere that it is possible to obtain the Lie algebroid
structure on $A\oplus B$ over $M$ directly from the two structures on
${\cal A}$.


\newcommand{\noopsort}[1]{} \newcommand{\singleletter}[1]{#1}

\end{document}